\documentclass{gtart_h}

\def\ifplaintex{\expandafter\ifx\csname documentclass\endcsname\relax}

\def\gtp{{\mathsurround=0pt\it $\cal G\mskip-2mu$eometry \&\ 
$\cal T\!\!$opology $\cal P\!$ublications}}  

\def\recd{{\small Received:\qua\receiveddate\ifx\reviseddate\relax
\else\qquad Revised:\qua\reviseddate\fi\par}} 

\def\lognumber#1{\def\thelognumber{#1}}
\def\volumenumber#1{\def\thevolumenumber{#1}}
\def\volumeyear#1{\def\thevolumeyear{#1}}
\def\papernumber#1{\def\thepapernumber{#1}}
\def\pagenumbers#1#2{\def\startpage{#1}\def\finishpage{#2}}
\def\published#1{\def\publishdate{#1}}

\def\received#1{\def\receiveddate{#1}}
\def\revised#1{\def\reviseddate{#1}}
\def\accepted#1{\def\accepteddate{#1}}

\def\asciiaddress#1{\def\theasciiaddress{#1}}
\def\asciiemail#1{\def\theasciiemail{#1}}

\long\def\asciiabstract#1{\long\def\theasciiabstract{#1}}

\let\\\par\let\thelognumber\relax\let\thevolumenumber\relax
\let\thepapernumber\relax\let\thevolumeyear\relax\let\startpage\relax
\let\finishpage\relax\let\publishdate\relax\let\receiveddate\relax
\let\reviseddate\relax\let\accepteddate\relax\let\theasciititle\relax
\let\theasciiauthors\relax\let\theasciiaddress\relax
\let\theasciiabstract\relax

\let\theasciiemail\relax

\ifplaintex
\font\logobig=cmssbx10 scaled 3836
\font\logomed=cmssbx10 scaled 2557
\else
\font\logobig=cmssbx10 scaled 4200
\font\logomed=cmssbx10 scaled 2800
\fi

\long\def\makeagttitle{   
\count0=\startpage
\agt\hfill      
\hbox to 45truept{\vbox to 0pt{\vglue -13truept{\logomed A\kern -.37em{\logobig 
T}\kern -.38em G}\vss}\hss}
\break
{\small Volume \thevolumenumber\ (\thevolumeyear)
\startpage--\finishpage\nl
Published: \publishdate}

\vglue .25truein

{\parskip=0pt\leftskip 0pt plus
1fil\def\\{\par\smallskip}{\Large\bf\thetitle}\par\medskip} \vglue
0.05truein

{\parskip=0pt\leftskip 0pt plus 1fil\def\\{\par}{\sc\theauthors}
\par\medskip}%
 
\vglue 0.03truein 

{\small\leftskip 25truept\rightskip 25truept{\bf Abstract}\stdspace\theabstract

{\bf AMS Classification}\stdspace\theprimaryclass
\ifx\thesecondaryclass\relax\else; \thesecondaryclass\fi\par
{\bf Keywords}\stdspace \thekeywords\par}\vglue 7truept

}   

\ifplaintex
\font\phead=cmsl9 scaled 950
\font\pnum=cmbx10 scaled 913
\font\pfoot=cmsl9 scaled 950
\headline{\vbox to 0pt{\vskip -4.5mm\line{\small\phead\ifnum
\count0=\startpage ISSN 1472-2739 (on-line) 1472-2747 (printed)
\hfill {\pnum\folio}\else\ifodd\count0\def\\{ }%
\ifx\theshorttitle\relax\thetitle\else\theshorttitle\fi\hfill{\pnum\folio}
\else\def\\{ and }{\pnum\folio}\hfill\ifx\theshortauthors\relax\theauthors
\else\theshortauthors\fi\fi\fi}\vss}}
\footline{\vbox to 0pt{\vglue 0mm\line{\small\pfoot\ifnum\count0=\startpage
\copyright\ \gtp\hfill\else
\agt, Volume \thevolumenumber\ (\thevolumeyear)\hfill\fi}\vss}}
\else
\font=cmsl9 scaled 1050
\font\lnum=cmbx10 
\font=cmsl9 scaled 1050
\makeatletter
\def\@oddhead{{\small\lhead\ifnum\count0=\startpage ISSN 1472-2739 
(on-line) 1472-2747 (printed)\hfill {\lnum\number\count0}\else\ifodd\count0
\def\\{ }\ifx\theshorttitle\relax \thetitle \else\theshorttitle\fi\hfill
{\lnum\number\count0}\else\def\\{ and }{\lnum\number\count0}
\hfill\ifx\theshortauthors\relax 
\theauthors\else\theshortauthors\fi\fi\fi}}\def\@evenhead{\@oddhead}
\def\@oddfoot{\small\lfoot\ifnum\count0=\startpage\copyright\ \gtp\hfill\else
\agt, Volume \thevolumenumber\ (\thevolumeyear)\hfill\fi}
\def\@evenfoot{\@oddfoot}
\makeatother
\fi
\let\maketitlepage\makeagttitle
\let\makeshorttitle\maketitlepage
\let\maketitle\maketitlepage

\newwrite\gtoutfile
\long\gdef\makeheadfile{  
{\def\\{, }\def\s{ }
\immediate\openout\gtoutfile head.xxx
\immediate\write\gtoutfile{Proxy-for: \ifx\theasciiauthors\relax
\theauthors\else\theasciiauthors\fi\s<\ifx\theasciiemail\relax\theemail\else\theasciiemail\fi>}
\immediate\write\gtoutfile{\noexpand\\}
\immediate\write\gtoutfile{Authors: \ifx\theasciiauthors\relax
\theauthors\else\theasciiauthors\fi}
{\def\\{ }\immediate\write\gtoutfile{Title: \ifx\theasciititle\relax
\thetitle\else\theasciititle\fi}}
\immediate\write\gtoutfile{Subj-class: GT or SG, GR etc}
\immediate\write\gtoutfile{MSC-class: \theprimaryclass\ifx\thesecondaryclass\relax\else, \thesecondaryclass\fi}
\immediate\write\gtoutfile{Journal-ref: Algebr. Geom. Topol. \thevolumenumber\s
(\thevolumeyear) \startpage-\finishpage}
\immediate\write\gtoutfile{Comments: Published by Algebraic and
Geometric Topology at}
\immediate\write\gtoutfile{\s\s\s  http://www.maths.warwick.ac.uk/agt/AGTVol\thevolumenumber/agt-\thevolumenumber-\thepapernumber.abs.html}
\immediate\write\gtoutfile{\noexpand\\}
\immediate\write\gtoutfile{}
\ifx\theasciiabstract\relax
\immediate\write\gtoutfile{\theabstract}\else
\immediate\write\gtoutfile{\theasciiabstract}\fi
\immediate\write\gtoutfile{}
\immediate\write\gtoutfile{\noexpand\\}
\immediate\write\gtoutfile{}
\immediate\closeout\gtoutfile}}  

\def\maketitlepage{\makeagttitle\makeheadfile}
\let\makeshorttitle\maketitlepage
\let\maketitle\maketitlepage

\lognumber{55}
\volumenumber{5}
\volumeyear{2005}
\papernumber{55}
\pagenumbers{1419}{1432}
\received{22 April 2004} 
\revised{13 September 2005}
\accepted{20 September 2005}
\published{15 October 2005}

\usepackage{graphicx,psfrag,amsmath,amssymb}

\def\psfraga <#1,#2> #3#4{%
\psfrag {#3}{\smash{\rlap{\kern #1 \raise #2\hbox{#4}}}}}

\def\figref#1{\hyperlink{#1anchor}{Figure~\ref*{#1}}}
\def\anchor#1{\noindent\hypertarget{#1anchor}{\smash{$\phantom{99}$}}\newline}
\newtheorem{theorem}{Theorem} [section]
\newtheorem{thm}[theorem]{Theorem}

\newtheorem{cor}[theorem]{Corollary}

\newtheorem{lemma}[theorem]{Lemma}

\newtheorem{prop}[theorem]{Proposition}

\begin{document}

\title{Intrinsically linked graphs and even linking number}
\author{Thomas Fleming\\Alexander Diesl}

\address{University of California San Diego, Department of 
Mathematics\\9500 Gilman Drive, La Jolla, CA 92093-0112, USA}
\secondaddress{University of California Berkeley, Department of 
Mathematics\\970 Evans Hall, Berkeley, CA 94720-3840, USA}
\asciiaddress{University of California San Diego, Department of 
Mathematics\\9500 Gilman Drive, La Jolla, CA 92093-0112, 
USA\\and\\University of California Berkeley, Department of 
Mathematics\\970 Evans Hall, Berkeley, CA 94720-3840, USA}
\asciiemail{tfleming@math.ucsd.edu, adiesl@math.berkeley.edu}
\gtemail{\mailto{tfleming@math.ucsd.edu}{\qua\rm and\qua}\mailto{adiesl@math.berkeley.edu}}

\begin{abstract}
We study intrinsically linked graphs where we require that every
embedding of the graph contains not just a non-split link, but a link
that satisfies some additional property. Examples of properties we
address in this paper are: a two component link with $lk(A,L) =
k2^{r}, k\neq0$, a non-split $n$-component link where all linking
numbers are even, or an $n$-component link with components $L, A_{i}$
where $lk(L,A_{i}) = 3k, k \neq 0$.  Links with other properties are
considered as well.
 
For a given property, we prove that every embedding of a certain
complete graph contains a link with that property.  The size of the
complete graph is determined by the property in question.
\end{abstract}

\asciiabstract{%
We study intrinsically linked graphs where we require that every
embedding of the graph contains not just a non-split link, but a link
that satisfies some additional property. Examples of properties we
address in this paper are: a two component link with lk(A,L) = k2^r, k
not 0, a non-split n-component link where all linking numbers are
even, or an n-component link with components L, A_i where lk(L,A_i) =
3k, k not 0.  Links with other properties are considered as well.  For
a given property, we prove that every embedding of a certain complete
graph contains a link with that property.  The size of the complete
graph is determined by the property in question.}

\primaryclass{57M15}
\secondaryclass{57M25,05C10}
\keywords{Intrinsically linked graph, spatial graph, graph embedding, 
linking number}

\maketitle

\section{Introduction}In the early 
1980's, Sachs, and independently, Conway and Gordon proved that every 
embedding of $K_{6}$ contains two disjoint cycles that form a nontrivial 
link \cite{sachs}, \cite{c&g}.  Since that time there has been extensive 
study of \emph{intrinsically linked graphs}, that is, graphs whose every 
embedding into three-space contains a non-split link.  The classification 
of these graphs, in terms of excluded minors, was settled by Robertson, 
Seymour and Thomas in \cite{r&s}.  Study of intrinsically linked graphs 
has turned to finding conditions that guarantee more complex structures 
within every embedding of the graph.  For example, one could require that 
every embedding of a graph contains two disjoint links \cite{foisyetal}, 
or a link with three or more components \cite{flapan2} \cite{foisyb}, or a 
link where one component is a nontrivial knot \cite{me}, or a two 
component link whose linking number is larger than some constant 
\cite{flapan3} \cite{s-t}.

It is this last case that interests us here.  We will prove that every embedding of certain complete graphs always contains links (of various numbers of components) whose linking numbers are multiples of some constant (generally a power of two). These results are Ramsey theoretic in nature, in that they imply certain structures arise in sufficiently large complete graphs.  In essence, we are striving to prove that given stringent requirements on a link, we can find a graph so that every embedding of that graph contains a link meeting those requirements.  

In the papers of Flapan \cite{flapan3} and Shirai and Taniyama \cite{s-t}, the authors study graphs that always contain two component links with linking number greater than or equal to $k$, where here, in our investigation of two component links, we demand that the linking number be a multiple of a power of two.  This is a much more restrictive condition and hence it should not be surprising that our bounds are much larger.

For two component links, we have the following theorems:

\medskip

\noindent \textbf{Theorem \ref{mod2whitehead}}\qua\textsl{Every embedding of 
$K_{\beta_{r}}$ contains a two component link with $lk(A,L)=2^{r}k$, for 
some $k 
\neq 0$.}

\medskip

\noindent \textbf{Theorem \ref{mod3thm}}\qua\textsl{Every embedding of $K_{35}$ 
contains a link of two components with $lk(A,L) = 3k$, for some $k \neq 
0$.}

\medskip

We can obtain a similar result for three component links, but we would like to see it strengthened so that all pairwise linking numbers are zero mod $2^{r}$.

\medskip

\noindent \textbf{Theorem \ref{mod2rings}}\qua\textsl{Every embedding of 
$K_{\delta_{r}}$ contains a non-split three component link $L,W,A$ with 
$lk(L,W) = 2^{r}k$ , $lk(L,A) = 2^{r}k'$ and $lk(W,A) \equiv 0$ mod $2$, 
for some $k,k' \neq 0$.} 

\medskip

There are also a number of results for links of arbitrarily many components.

\medskip

\noindent \textbf{Corollary \ref{mod2keys}}\qua\textsl{Every embedding of 
$K_{\beta'_{n,r}}$ contains an $n+1$ component link \\ $L, Z_{1}, \ldots, 
Z_{n}$ with $lk(L,Z_{i}) = 2^{r}k$, for some $k \neq 0$.}

\medskip

\noindent \textbf{Corollary \ref{mod3keys}}\qua\textsl{Every embedding of 
$K_{7\alpha_{3n}'}$ contains an $n+1$ component link \\ $L, Z_{1}, \ldots, 
Z_{n}$ with $lk(L,Z_{i}) = 3k$, for some $k \neq 0$.}

\medskip

\noindent \textbf{Theorem \ref{alleven}}\qua\textsl{Every embedding of $K_{\epsilon_{n}}$ contains a non-split link of n+1 components where all the pairwise linking numbers are even.}

\medskip

The sequence $\gamma_{n}= \Pi_{i=0}^{n-1}(2^{i}+1)$ plays a significant role in determining the size of a complete graph required for some of these properties.  Unfortunately, $\gamma_{n}$ grows faster than $n^{n}$.  Obviously, we would like to see better bounds on the number of vertices required.  We also note that in this paper we provide upper bounds for the number of vertices required for these behaviors.  Obvious lower bounds are seven for all cases \cite{c&g}, and ten for cases involving links of three or more components \cite{flapan}.  Finding better lower bounds or even the minimum number of vertices required for these properties remain open questions.

\section{The basic linking construction}

Our primary 
technique in this paper is to take two cycles that link with a third and 
use them to produce new cycles with various properties.  To that end, the 
following technical lemmas will be used extensively in later results.

\begin{lemma}Let a graph $G$ be embedded such that there exists a 
three component link with components $L, Z, W$ which have 
$lk(L,Z)=q_{1}>0$, $lk(L,W)=q_{2}>0$ with $q_{i} \equiv 0$ mod $2^{r}$. 
Assume that there exist $2^{r+1}+1$ paths $P_{i}$ from $Z$ to $W$ with 
distinct end points such that the interior of the $P_{i}$ are disjoint 
from $L, W, Z$ and from the other $P_{i}$. Further assume that when 
traversing $W$ in the direction of its orientation, we meet the $P_{i}$ in 
ascending order, and when traversing $Z$ in a similar fashion, we meet 
them in descending order. Then this embedding of $G$ contains a two 
component link $A \cup L$ with $lk(A,L)=k 2^{r+1}$, for some $k \neq 0$.  
\label{evenlinklemma}
\end{lemma}
\begin{proof}
Clearly, if $q_{1}$ 
or $q_{2}$ is zero mod $2^{r+1}$, we are done, so we may assume that 
$q_{1} \equiv q_{2} \equiv 2^{r}$ mod $2^{r+1}$.Orient the $P_{i}$ 
such that going from $Z$ to $W$ is the positive direction.Let $A_{1}$ be 
the oriented simple closed curve formed by traversing $P_{1}$ in the 
positive direction, along $W$ in the positive direction from $P_{1}$ to 
$P_{2}$, in the negative direction along $P_{2}$, and finally along $Z$ in 
the positive direction from $P_{2}$ to $P_{1}$.  See 
\figref{aloops}. Define the other $A_{i}$ similarly.

\begin{figure}[ht!]\small\anchor{aloops}
\psfrag {1}{$1$}
\psfrag {2}{$2$}
\psfrag {3}{$3$}
\psfrag {Z}{$Z$}
\psfrag {W}{$W$}
\psfrag {2Z}{$2_Z$}
\psfrag {1Z}{$1_Z$}
\psfrag {1W}{$1_W$}
\psfrag {2W}{$2_W$}
\psfrag {A1}{$A_1$}
\cl{\includegraphics[width=2.5in]{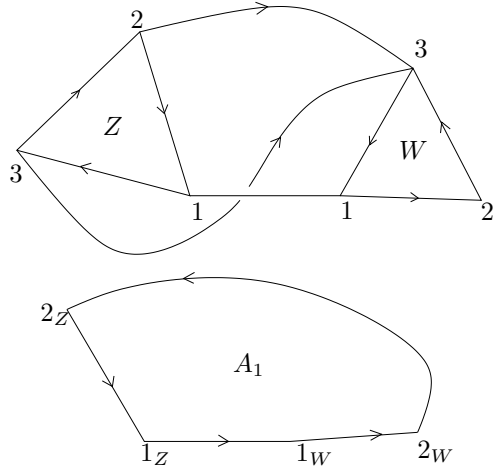}}
\caption{Constructing the $A_{i}$ for r=1}
\label{aloops}\end{figure}

Note that now, we 
have that in $H_{1}(\mathbf{R}^{3}-L;\mathbf{Z})$;$$[Z] + [W] = 
\sum [A_{i}]$$Notice $[Z] + [W] = q_{1} + q_{2}$, which is 0 mod 
$2^{r+1}$ and non-zero integrally.

For $1 \leq i \leq 2^{r+1}$, let $S_{i} = \sum_{j=1}^{i} [A_{i}]$.  

Suppose $S_{i} \equiv 0$ mod $2^{r+1}$ for some $i$. If $S_{i} \neq 0$ integrally, let $A = A_{1} \cup A_{2} \cup \ldots \cup A_{i}$.  If $S_{i} = 0$ integrally, then $\sum_{i+1}^{2^{r+1}+1} [A_{j}] = q_{1} + q_{2}$. Thus by setting $A = A_{i+1} \cup \ldots \cup A_{2^{r+1}+1}$, we have the desired result.

Thus we may now assume that $S_{i} \neq 0$ mod $2^{r+1}$ for all $i$.  
However, we have $2^{r+1}$ values, so for some $i < j$, we have that 
$S_{i} \equiv S_{j}$ mod $2^{r+1}$.  Thus $0 \equiv S_{j} - S_{i} \equiv 
\sum_{k=i+1}^{j} [A_{k}]$ mod $2^{r+1}$.  Now, as before, if $S_{j} - 
S_{i} \neq 0$ integrally, let $A = A_{i+1} \cup \ldots A_{j}$.  If $S_{j} 
- S_{i} = 0$ integrally, then let $A = A_{1} \cup \ldots \cup A_{i} \cup 
A_{j+1} \cup \ldots \cup A_{2^{r+1}+1}$.
\end{proof}

Define $\alpha_{1} = 6$, $\alpha_{2} = 10$, and for $m \geq 1, 
\alpha_{2m+1} = 2\alpha_{2m} + 6$ and $\alpha_{2m+2} = 2\alpha_{2m+1}$.  
Then for $m \geq 1$ we have $\alpha_{2m+1} = 6(\sum_{j=0}^{m} 4^{j}) - 
4^{m}$.  
Note that $\alpha_{1} = 6\alpha_{1}'$ and for $n \geq 2$, $\alpha_{n} < 6\alpha_{n}'$, where $\alpha_{n}'$ is the sequence defined before Corollary \ref{starcor}.

\begin{lemma}Every embedding of $K_{\alpha_{n}}$ contains a nonsplit link of n+1 components labeled $L, Z_{1}, \ldots, Z_{n}$, where $lk(L,Z_{i}) \neq 0$ for all $i$.  
\label{ringofkeys}
\end{lemma}

\begin{proof}
The case $n=1$ was shown in \cite{c&g}.

We will proceed by induction.  The base case is $n = 2$, and the desired result for $K_{10}$ was shown in \cite{flapan}. The inductive step for the even and odd cases must be handled separately.

Now suppose that $K_{\alpha_{n-1}}$ contains such of link of $n$ components, where $n$ is even. The graph $K_{\alpha_{n}}$ is $K_{2\alpha_{n-1}}$ and we may partition it into two disjoint copies of $K_{\alpha_{n-1}}$.  Thus we have found two cycles, $L$ and $L'$ which each have non-zero linking number with $n-1$ other cycles. Clearly if $lk(Z_{i},L') \neq 0$ or $lk(Z_{i}',L)\neq 0$ for some $i$, we are done.  So we may assume that all such linking numbers are zero. 

Now, choose two vertices $l_{1}$ and $l_{2}$ in $L$ and two vertices $l_{1}'$ and $l_{2}'$ in $L'$. Let $P$ be a path in $L$ connecting $l_{1}$ to $l_{2}$, and $P'$ be a path in $L'$ connecting $l_{1}'$ to $l_{2}'$. Since this is a complete graph, there are edges connecting $l_{1}$ to $l_{1}'$ and $l_{2}$ to $l_{2}'$.  Note that the interiors of these edges are disjoint from $L, L', Z_{i},$ and $Z_{i}'$.   Let $V$ denote the cycle formed by these edges and $P$ and $P'$.  See \figref{Vpic}. 
\begin{figure}[ht!]\small\anchor{Vpic}
\psfrag {P}{$P$}
\psfrag {V}{$V$}
\psfrag {L}{$L$}
\psfrag {Pp}{$P'$}
\psfrag {Lp}{$L'$}
\psfrag {l1}{$l_1$}
\psfrag {l2}{$l_2$}
\psfrag {l1p}{$l'_1$}
\psfrag {l2p}{$l'_2$}
\cl{\includegraphics[width=2.5 in]{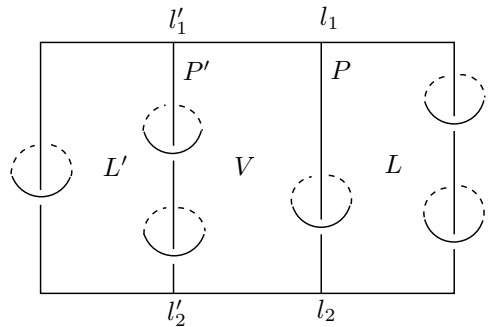}}
\caption{The cycles $L$, $L'$ and V in the case $n=4$.  The $Z_{i}$ are shown as ovals.}
\label{Vpic}
\end{figure}

Now, examine $lk(V,Z_{i})$ and $lk(V,Z_{i}')$. If $n$ of these are 
non-zero, then we have constructed the desired link.  If fewer than $n-1$ of these are non-zero, then we take the cycle $L \cup V \cup L'$, and since linking number is additive, we have the desired link once again.

If exactly $n-1$ of the $Z_{i}, Z_{i}'$ have non-zero linking number with 
$V$, then as $n-1$ is odd, more than half of them came from one set, say 
the $Z_{i}'$.  Again using the additivity of linking number, the cycle $V 
\cup L$, has $lk(V \cup L, Z_{i}) \neq 0$ for at least $\frac {n}{2}$ of 
the $Z_{i}$, and similarly for the $Z_{i}'$.  This gives a link with $n+1$ 
components with the desired properties when $n$ is even.

When $n$ is odd, the inductive step requires more care. In the graph $K_{\alpha_{n}}$ we may find a two copies of $K_{\alpha_{n-1}}$ as well as a copy of $K_{6}$.  The latter will tip the balance in our favor.

Again choose vertices $l_{i}$ and $l_{i}'$.  Choose also vertices $x$ and $y$ on one of the two linked triangles in $K_{6}$, call this triangle $L''$, the other $T$.  Let $P$ and $P'$ be the paths in $L$ and $L'$ as before.  Now let $V$ be the cycle formed by $P$, the edge $l_{2}-l_{2}'$, $P'$, $l_{1}'-x$, $x-y$, $y-l_{1}$.  See \figref{Vfig2}.

\begin{figure}[ht!]\small\anchor{Vfig2}
\psfrag {V}{$V$}
\psfrag {T}{$T$}
\psfrag {L}{$L$}
\psfrag {Lp}{$L'$}
\psfrag {Ld}{$L''$}
\cl{\includegraphics[width=2.5 in]{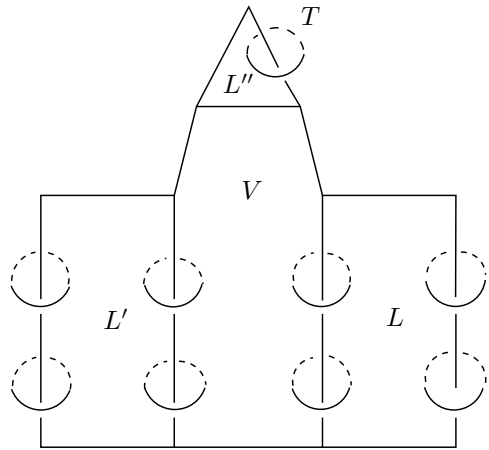}}
\caption{The cycles $L$, $L'$, $L''$ and V in the case $n=5$.  The $Z_{i}$ and $T$ are shown as ovals.}
\label{Vfig2}
\end{figure}

Now, clearly $lk(Z_{i}, L')$, $lk(Z_{i}', L)$, $lk(L,T)$, and $lk(L',T)$ must be zero, or we have the desired link. Once again, if $n$ or more of the $Z_{i}, Z_{i}'$ have non-zero linking number with $V$ we are done, and if fewer than $n-1$, then we take $L \cup V \cup L'$ as before.  This leaves only the case where $n-1$ of the $Z_{i}, Z_{i}'$ have non-zero linking number with $V$.  

If $lk(T,V)\neq 0$, then take the curve $V$.  This has non-zero linking number with $n-1$ of the $Z_{i}, Z'_{i}$, and non-zero linking with $T$, so we are done in this case as well. 

Thus, we may assume $lk(T,V)$ is zero. Let $C$ be the subset of $\{Z_{i}, Z_{i}'\}$ which have zero linking number with $V$, and $D$ be the subset which have non-zero linking number with $V$.  Suppose $r$ elements of $C$ and $s$ elements of $D$ have non-zero linking with $L''$.  Of these $s$ elements, suppose $t$ satisfy the equation $lk(Z_{i},L) + lk(Z_{i},V) + lk(Z_{i}, L'') = 0$ (or the similar version for $Z_{i}', L'$).  Then at least $n + r - s + t$ elements have non-zero linking with $V \cup L''$ and at least $n - r + s - t$ have non-zero linking with $V \cup L \cup L' \cup L''$.  Clearly one of these is greater than or equal to $n$, so we are done.
\end{proof}

\begin{lemma}Every embedding of $K_{\alpha_{n}}$ contains a nonsplit link of n+1 components labeled $L, Z_{1}, \ldots, Z_{n}$, where $lk(L,Z_{i}) \equiv 1$ mod $2$ for all $i$.  
\label{ringofkeys2}
\end{lemma}
The proof of Lemma \ref{ringofkeys2} proceeds exactly as for Lemma \ref{ringofkeys}, replacing the phrase (non-)zero with (non-)zero mod two.   

Let $\alpha_{1}' =1$, $\alpha_{2m}' = 2\alpha_{2m-1}'$, and $\alpha_{2m+1}'=2\alpha_{2m}'+1$.  Again, a calculation gives $\alpha_{2m-1}' = \frac{4^{m}-1}{3}$.

\begin{cor}
If every embedding of $G$ contains a two component link with non-zero 
linking number, then every embedding of $*^{\alpha_{n}'}G$ contains an $n+1$ component link $L, Z_{i}$ with $lk(L,Z_{i}) \neq 0$. 
\label{starcor}
\end{cor}

\begin{proof}
Note that the sequence $\alpha_{n}'$ satisfies the following recursive relations: 
$\alpha_{n}' = 2\alpha_{n-1}'$ if $n$ is even, $\alpha_{n}' = 2\alpha_{n-1}' + 1$ if $n$ is odd, and $\alpha_{1}'=1$.

The proof is now the same as the proof of Lemma \ref{ringofkeys}, but replace $K_{10}$ with $G*G$ and $K_{6}$ with $G$.
\end{proof}

Naturally, for Lemma \ref{ringofkeys2} we have a corollary similar to Corollary \ref{starcor}.

We will now obtain a sequence of links that can be thought of as the ``mod two Whitehead links."

Let $\beta_{0} = 6$, $\beta_{1} = 10$ and $\beta_{r} = \alpha_{2^{r}}'(\gamma_{r}+3)$ for higher $r$. Here $\gamma_{r} = \Pi_{i=0}^{r-1}(2^{i}+1)$.

\begin{thm}Every embedding of $K_{\beta_{r}}$ contains a two component 
link with $lk(A,L)=2^{r}k$, for some $k \neq 0$. \label{mod2whitehead}
\end{thm} 

\begin{proof}

Notice that the case $r=0$ was proved for $K_{6}$ in \cite{c&g}.

The case $r = 1$ is direct. 
By \cite{flapan} every embedding of $K_{10}$ contains a three component link, with $lk(L,Z)=q_{1}, lk(L,W)=q_{2}$ with $q_{i}$ odd. Choose an orientation of $L \cup W \cup Z$ so that the $q_{i}$ are positive. Clearly $Z$ and $W$ contain at least three vertices. Since $K_{10}$ is a complete graph, we have edges from each vertex of $Z$ to all vertices of $W$, and the interiors of these edges are disjoint from $L,W,Z$.  We may choose $P_{i}$ from among these edges that satisfy the conditions of Lemma \ref{evenlinklemma}, so we can find the desired link.

We will now give the construction for general $r$.  

Every embedding of $K_{m+3}$ contains a link composed of a triangle and an $m$-cycle with nonzero linking number by \cite{johnson}.  By Corollary \ref{starcor}, we know that $K_{m+3}*K_{m+3}$ contains a three component link where at least two of the three pairwise linking numbers are nonzero.  Choosing the triangles to act as $L$, we know that the $Z_{i}$ have $m$ vertices. 

Thus, in $K_{2m+6}$ we may find a three component link $L, Z_{1}, Z_{2}$ with $lk(L,Z_{i}) \neq 0$, and with each $Z_{i}$ containing $m$ vertices.  Now, by Corollary \ref{starcor} and the preceding discussion, if we want a $n+1$ component link with these properties, we need only take $K_{\alpha_{n}'(m+3)}$.

Suppose $n$ is a large power of two, and $m$ is large.  Pair up the 
$Z_{i}$. If $lk(Z_{i},L)$ is even for some $i$, ignore that pair.  For the 
remaining pairs (where $lk(Z_{i},L) \equiv 1$ mod 2)  choose three evenly 
spaced paths between the members of that pair.  Applying Lemma 
\ref{evenlinklemma} to these pairs, we have new cycles $A_{i}$ with 
$lk(L,A_{i})=2k$ ($k \neq 0$).  Now including the $Z_{i}$ with even 
linking number, we have $\frac{n}{2}$ cycles which have nonzero even 
linking number with $L$ and each cycle has at least $2(\frac{m}{3}+1)$ 
vertices. As long as we have chosen $m$ and $n$ large enough, we may 
iterate 
this process to obtain the desired link.That is, suppose $lk(Z_{1},L) 
\equiv lk(Z_{2},L) \equiv 0$ mod $2^{s-1}$ but are non-zero integrally, 
and that the $Z_{i}$ have at least $m'$ vertices.  If $lk(Z_{i},L) \equiv 0$ 
mod $2^{s}$, we are done.  So we may assume that $lk(Z_{1},L) \equiv 
lk(Z_{2},L) \equiv 2^{s-1}$ mod $2^{s}$.  Choosing $2^{s}+1$ paths from 
$Z_{1}$ to $Z_{2}$ allows us to apply Lemma \ref{evenlinklemma} to obtain 
a cycle $Z$ with $lk(Z,L) = k2^{s}, k \neq 0$, such that $Z$ has at least 
$2(\frac{m'}{2^{s}+1}+1)$ vertices.

Thus, if we want a two component link with $lk(A,L) \equiv 0$ mod $2^{r}$, 
we choose $n = 2^{r}$, and $m = \gamma_{r}$.  Clearly we need $n=2^{r}$ 
components linked with $L$, as after each pairing, we have increased the 
linking number of $Z_{i}$ with $L$ by a power of two, but have halved the 
number of components.  However, on iteration $s$ in this procedure (moving 
from $2^{s-1}$ to $2^{s}$) we need to have at least $2^{s}+1$ vertices in 
each $Z_{i}$ to successfully apply Lemma \ref{evenlinklemma}.  Since the 
number of vertices is strictly decreasing in this process, it suffices to 
guarantee that we have $c_{r}$ vertices at the $r$th iteration, where 
$c_{r} \geq 2^{r} +1$.  If we have $c_{s}$ vertices at the beginning of 
iteration $s$, we will be left with $c_{s+1} = 
2(\frac{c_{s}}{2^{s}+1}+1)$ vertices at the beginning of iteration $s+1$. Thus, when $c_{s+1}$ is even, we have $c_{s} = 
(2^{s}+1)(\frac{c_{s+1}}{2}-1)$, and we want to choose $m \geq c_{1}$.  
Choosing $c_{r} = 2^{r} + 2$ gives $c_{r-1} = (2^{r-1}+1)(2^{r-1})$.  
Since we are looking only for a bound, let us take $c_{r-1} = 
(2^{r-1}+1)(2^{r-1}) + 2$.  This gives $c_{r-2} = 
(2^{r-2}+1)(2^{r-1}+1)(2^{r-2})$.  Continuing this process gives $c_{1} = 
\gamma_{r}$.  Thus, by taking $m = \gamma_{r}$ we ensure that at each 
stage, we have a sufficient number of vertices to apply Lemma 
\ref{evenlinklemma}.
\end{proof}

We have calculated the actual number of vertices required for some small cases, and $\gamma_{r}$ provides a very reasonable estimate.  Denoting the required number of vertices $v_{r}$, we note that $\gamma_{1} = v_{1} = 3$; $\gamma_{2} = v_{2} = 6$;$\gamma_{3} = 30, v_{3} = 27$;$\gamma_{4} = 270, v_{4} =261$.

Let $\beta_{n,r}' = \alpha_{n2^{r}}'(\gamma_{r}+3).$

\begin{cor} Every embedding of $K_{\beta'_{n,r}}$ contains an $n+1$ 
component link \\ $L, Z_{1}, \ldots, Z_{n}$ with $lk(L,Z_{i}) = 2^{r}k$, 
for 
some $k \neq 0$. \label{mod2keys} \end{cor}

\begin{proof}

By Lemma \ref{ringofkeys}, we can find an $n2^{r}+1$ component link where 
$lk(Z'_{i},L) \neq 0$ and each of the $Z_{i}$ has at least $\gamma_{r}$ vertices.  Applying the construction of Theorem \ref{mod2whitehead} to each set of $2^{r}$ components, we obtain the desired link.  Note that the construction of one such $Z_{i}$ does not affect $L$ or the other $Z_{i}$ so we may treat them independently.
\end{proof}

Note that Theorem \ref{mod2whitehead} implies that every embedding of $K_{10}$ contains a non-split two component link with even linking number.  In fact, we can improve upon that result.  

\begin{prop}
Every embedding of $K_{10}$ contains a two component link $A, L$ with $lk(A,L) \equiv 2$ mod $4$. 
\label{k10prop}
\end{prop}

\begin{proof}
The proof of Thoerem 2.4 in \cite{menp} shows that every embedding of $K_{10}$ contains a two component link with odd unoriented Sato-Levine invariant.  A pair of intersection points of opposite sign contribute an even quantity to $\beta^{*}$; a pair of intersection points with matching sign contribute an odd amount.  Thus, it is easy to see that $\beta^{*}(A,L)$ is even if and only if $lk(A,L) \equiv 0$ mod 4.   Since the constructed link has odd unoriented Sato-Levine invariant, it must have linking number 2 mod 4.  
\end{proof}

\section{Multiple component constructions}

We would like to see the following theorem extended to have all three pairwise linking numbers be zero mod $2^{r}$.

Let $\delta_{r} = \alpha'_{3(2^{r})}((2^{2r-1}+2^{r})\gamma_{r} + 3)$. 

\begin{thm}Every embedding of $K_{\delta_{r}}$ contains a non-split three 
component link $L,W,A$ with $lk(L,W) = 2^{r}k$ , $lk(L,A) = 2^{r}k'$ and 
$lk(W,A) \equiv 0$ mod $2$, for some $k,k' \neq 0$. \label{mod2rings}
\end{thm} 

\begin{proof}

By Corollary \ref{mod2keys}, we may find a four component link $L,W,Z_{1}, Z_{2}$, with $lk(L,W) \equiv lk(L,Z_{i}) \equiv 0$ mod $2^{r}$, but nonzero integrally.  In building up the $Z_{i}$ as we did in the proof of Theorem \ref{mod2whitehead}, at each stage the linking number is multiplied by two, but the number of vertices is reduced from $c_{j}$ to $c_{j+1}$.  If we had begun with $\gamma_{r}$ vertices in each cycle, the construction would end with the $Z_{i}$ containing at least four vertices.  However, we will require $(2^{r}+1)^{2}$ vertices, so we multiply $\gamma_{r}$ by the smallest constant that ensures this result.
Thus, each of $W, Z_{i}$ has at least $2(2^{2r-1}+2^{r}+1) > (2^{r}+1)^{2}$ vertices. This number of vertices is necessary to form a cycle $A$ from the $Z_{i}$ so that the three component link $L,W,A$ will have the desired property. This number of vertices is not necessary in $W$, but unfortunately Lemma \ref{ringofkeys} does not allow us to control the size of the components independently.

If $lk(Z_{i},W) \equiv 0$ mod $2$ we are done, so suppose that $lk(Z_{i},W) \equiv 1$ mod $2$.  The $Z_{i}$ have $(2^{r}+1)^{2}$ vertices by construction.  Find $(2^{r}+1)^{2}$ paths $P_{i}$ and construct cycles $A_{i}'$ as in the proof of Lemma \ref{evenlinklemma}.  Now, write these in a list.  Group the first $2^{r}$ cycles, then leave one, then group the next $2^{r}$ cycles, etc.
Look at the first $2^{r}$ cycles in $H_{1}(\mathbf{R}^{3}-L;\mathbf{Z})$; we know that there is some subsequence whose sum is 0 mod $2^{r}$.  Glue these cycles together.  Look at the next block of $2^{r}$ cycles.  Find its subsequence.  Now identify all cycles between these two.  The result will be a sequence of cycles $A_{i}'$ where $[A_{i}'] \equiv 0$ mod $2^{r}$ for $i$ odd, and takes arbitrary values for $i$ even.  Any cycles coming before $A_{1}'$ may be pushed around to the end of the sequence since it is cyclic.

The result is a list of $2(2^{r}+1)$ cycles, half of which are zero mod $2^{r}$.  Now, ignoring the zeros, there is a list of $2^{r}+1$ numbers, and as before there must be a proper consecutive subsequence that is zero mod $2^{r}$.  Let $A''_{2}$ be the cycle which corresponds to this sequence, $A_{1}''$ be the cycle which is 0 mod $2^{r}$ that immediately precedes it, and $A_{3}''$ be the complement of these two cycles.  Since $\sum [A_{i}'']_{L} = [Z_{1}]_{L} + [Z_{2}]_{L} = k2^{r}, k \neq 0$, we know $[A_{3}'']_{L} \equiv 0$ mod $2^{r}$.

Now examine $[A_{i}''] \in H_{1}(\mathbf{R}^{3}-W;\mathbf{Z})$.  Again, there is a subsequence of these that sums to zero mod 2. Call the corresponding cycle $A$. Since $\sum [A_{i}'']_{W} = [Z_{1}]_{W} + [Z_{2}]_{W} \equiv 0$ mod 2, the complement of $A$ is also zero mod 2. Call it $A'$. 

By construction, $lk(A,L) \equiv 0$ mod $2^{r}$ and $lk(A,W) \equiv lk(A',W) \equiv 0$ mod 2.  If $lk(A,L) = 0$, then as $[A]_{L} + [A']_{L} = \sum [A_{i}'']_{L}$, we have that $lk(A',L) = k2^{r}, k \neq 0$ and we are done.      
\end{proof}

Let $\epsilon_{1} = 10$ and $\epsilon_{n} = \alpha'_{2^{n+1}-2}(\gamma_{n}' + 3)$. Here $\gamma_{n}' = \Pi_{i=1}^{n} 3(2^{i-1}) = 3^{n}2^{\frac{n(n-1)}{2}}$.

Then extending our construction yet again, we can show the following.

\begin{thm}Every embedding of $K_{\epsilon_{n}}$ contains a non-split link of n+1 components where all the pairwise linking numbers are even.  
\label{alleven}
\end{thm}

\begin{proof}

The case $n=1$ is proved in Lemma \ref{evenlinklemma} and $n=2$ is proved in Theorem \ref{mod2rings}.
For larger $n$, we may think of $K_{\epsilon_{n}}$ as the star of $\alpha_{2^{n+1}-2}'$ copies of $K_{\gamma_{n}'+3}$, and so it contains a link of $2^{n+1}-1$ components, labeled $L, Z_{i}$, where $lk(L, Z_{i}) \neq 0$ and the $Z_{i}$ have $\gamma_{n}'$ vertices. Since $2^{n+1}-2 = \sum_{i=1}^{n}2^{i}$, we may break the $Z_{i}$ into sets of size $2^{i}$ for each $1 \leq i \leq n$.

Take $Z_{1}$ and $Z_{2}$.  If $lk(Z_{i},L) \equiv 0$ mod 2, label that cycle $V_{1}$ and move on.  If $lk(Z_{i},L) \equiv 1$ mod 2 for $i=1,2$, choosing three paths between them and forming $A_{i}$ as usual, we may produce $V_{1}$, which has $lk(V_{1},L) \equiv 2k_{1}$, with $k_{1} \neq 0$.  

Take $Z_{3}, Z_{4}, Z_{5}, Z_{6}$.  Pair these and construct $Z_{3}', 
Z_{4}'$ 
which 
have even (and positive) linking number with $L$ as we just did for $Z_{1}$ and $Z_{2}$.  If $lk(Z_{i}',V_{1}) \equiv 0$ mod 2, label it $V_{2}$ and continue to the next set of $Z_{i}$.   If not, choose six paths from $Z_{3}'$ to $Z_{4}'$ and construct $A_{1}, \dots, A_{6}$.  Examine these cycles in $H_{1}(\mathbf{R}^{3}-L;\mathbf{Z/2})$.  Since $lk(Z_{i}',L) \equiv 0$ mod 2 this is a binary sequence which sums to zero.  Cut this sequence into the maximum number of blocks, where the sum in each block is zero.  Since we can cyclicly permute the entries, there must be at least three such blocks.  Produce the cycles $A_{1}', A_{2}', A_{3}'$ by fusing the cycles in each block.  Look at the $A_{s}'$ in $H_{1}(\mathbf{R}^{3}-V_{1};\mathbf{Z/2})$.  There is a proper subsequence whose sum is zero mod two.  Call the union of those cycles $A$, and their complement $A'$.  Now, $lk(A,V_{1}) \equiv lk(A',V_{1}) \equiv 0$ mod 2.  If $lk(A,L) = 0$, then as $[A]_{L} + [A']_{L} = [Z_{3}'] + [Z_{4}']$, we have that $lk(A',L) = 2k_{2}, k_{2} \neq 0$.  
Call this cycle $V_{2}$.

Take $Z_{(\sum_{s=1}^{j-1}2^{s})+1} \ldots Z_{\sum_{s=1}^{j}2^{s}}$. Pair 
them to produce $Z^{(1)}_{i}$ that have even, non-zero linking number with 
$L$ as we did for $Z_{1}$ and $Z_{2}$.  Now, pair the $Z^{(1)}_{i}$ and 
produce $Z^{(2)}_{i}$ that have $lk(Z_{i}^{(2)},V_{1}) \equiv 0$ mod 2, 
and $lk(Z_{i}^{(2)},L) = k2, k \neq 0$ as we did for $Z_{3}, Z_{4}, Z_{5}, 
Z_{6}$.  Now, pair the $Z^{(2)}_{i}$.  If either member of a pair has 
$lk(Z_{i}^{(2)},V_{2}) \equiv 0$ mod 2, it advances to the next round.  
Otherwise, choose twelve paths connecting the members of each pair and 
produce $A_{1} \ldots A_{12}$.  Examine these in 
$H_{1}(\mathbf{R}^{3}-L;\mathbf{Z/2})$, then 
$H_{1}(\mathbf{R}^{3}-V_{1};\mathbf{Z/2})$, breaking it into the maximum 
number of blocks each time.  We will be left with at least three cycles 
that we examine in $H_{1}(\mathbf{R}^{3}-V_{2};\mathbf{Z/2})$, allowing us 
to construct $Z^{(3)}_{i}$.

Continue the process, connecting the paired cycles by $3(2^{k-1})$ paths 
in the $k$th step and producing the $A_{s}$.  We have chosen our cycles 
$Z_{i}$ to have $\gamma_{n}'$ vertices, which ensures that we may find 
these paths in up to $n$ iterative steps. Examine these $A_{s}$ in 
$H_{1}(\mathbf{R}^{3}-L;\mathbf{Z/2})$, then 
$H_{1}(\mathbf{R}^{3}-V_{i};\mathbf{Z/2})$ for $i < k-1$ fusing cycles as 
necessary each time.  We began with $3(2^{k-1})$ cycles labeled $A_{s}$, 
so there must be at least three distinct cycles left after fusing $k-1$ 
times.  Examine these three cycles in 
$H_{1}(\mathbf{R}^{3}-V_{k-1};\mathbf{Z/2})$.  There is a subsequence that 
sums to zero, we take the corresponding cycle $W_{1}$ and its complement 
$W_{2}$.  Now, by construction, $lk(W_{l},V_{i}) \equiv 0$ mod 2 for $i < 
k$, and at least one of $lk(W_{l},L)$ is even and nonzero.  Choose that 
$W_{l}$ to be $Z_{i}^{(k)}$.

Since we began with $2^{j}$ cycles, we may repeat this process $j$ times until we have $V_{j}$, a cycle with $lk(V_{j},V_{i}) \equiv 0$ mod 2 for $i<j$ and $lk(V_{j},L) = 2k_{j}$ with $k_{j} \neq 0$.  

Iterating this process produces the desired link. 
\end{proof}

\section{Mod 3 constructions}

It is unclear if similar constructions will work modulo integers other 
than 2.  We provide a first step in that direction.

\begin{thm}Every embedding of $K_{35}$ contains a link of two components with $lk(A,L) = 3k$, $k \neq 0$.
\label{mod3thm}
\end{thm}

\begin{proof}

The graph $K_{35}$ has $*^{5}K_{7}$ as a subgraph. By Corollary \ref{starcor} and the discussion in the proof of Theorem \ref{mod2whitehead}, every embedding of $K_{35}$ contains a four component link $L, Z_{1}, Z_{2}, Z_{3}$ with $lk(L,Z_{i}) \neq 0$, and the $Z_{i}$ each have four vertices.  Orient the $Z_{i}$ and $L$ so that $lk(L,Z_{i}) = q_{i}> 0$.

If $lk(L,Z_{i}) \equiv 0$ mod $3$ for some $i$, then we are done.

If $lk(L,Z_{i}) \equiv 1$ mod $3$ and $lk(L,Z_{j}) \equiv 2$ mod $3$, then we may construct four cycles $A_{i}$ as in Lemma \ref{evenlinklemma}. Here $q_{1} + q_{2} \equiv 0$ mod $3$, and taking $S_{1}, S_{2}, S_{3}$ as before, the construction of $A$ proceeds in exactly the same manner.

Thus, we may assume that $lk(L,Z_{1}) \equiv lk(L,Z_{2}) \equiv lk(L,Z_{3})$ mod $3$.

We now construct cycles $A_{1}, A_{2}, A_{3}, A_{4}, A_{5}$ as in \figref{squaresfig}.  Note that in $H_{1}(\mathbf{R}^{3}-L;\mathbf{Z})$, 
$$\sum [Z_{i}] = \sum [A_{i}] = q_{1} + q_{2} + q_{3} \equiv 0 \quad 
\mbox{mod 3}$$

\begin{figure}[ht!]\small\anchor{squaresfig}
\psfrag {Z1}{$Z_1$}
\psfrag {Z2}{$Z_2$}
\psfrag {Z3}{$Z_3$}
\psfrag {A1}{$A_1$}
\psfrag {A2}{$A_2$}
\psfrag {A3}{$A_3$}
\psfrag {A4}{$A_4$}
\psfrag {A5}{$A_5$}
\cl{\includegraphics[width=2 in]{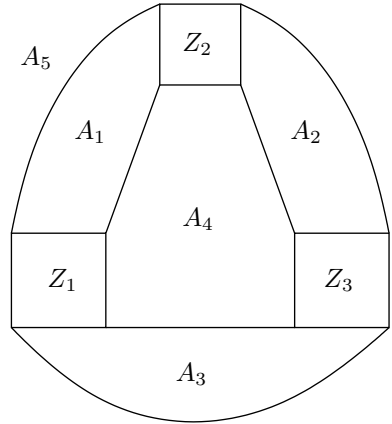}}
\caption{The cycles $A_{i}$ when using all the $Z_{i}$. The $A_{i}$ (except $A_{5}$) are oriented counter-clockwise, $A_{5}$ and the $Z_{i}$ clockwise.}
\label{squaresfig}
\end{figure}

If $[A_{i}] \equiv 0$ for some $i \geq 4$, then we take $A = A_{i}$ if $[A_{i}] \neq 0$.  If  $[A_{i}] = 0$ integrally, then let $A = \cup_{j \neq i} A_{j}$.

Suppose that $[A_{i}] \equiv 0$ for some $i \leq 3$. Again, take $A = A_{i}$ if $[A_{i}] \neq 0$ integrally.  If $[A_{i}] = 0$, look at the remaining $A_{j}$.  We have four contiguous cycles whose sum is zero mod 3, see \figref{squaresfig}. Forming the partial sums $S_{1}, S_{2}$, and $S_{3}$, either two of these are equal, or one is zero, and we may proceed as in Lemma \ref{evenlinklemma}.

Now, if $[A_{i}] \neq 0$ mod $3$ for all $i$, then two of $[A_{1}], 
[A_{2}], [A_{3}]$ are equal.  Say $[A_{1}] = [A_{2}]$. Take the following 
sums.  Let $S_{1} = [A_{4}]$, $S_{2} = [A_{4}] + [A_{1}]$, $S_{3} = 
[A_{4}] + [A_{1}] + [A_{5}]$, all taken mod 3.  We can ignore the cases 
$S_{1} \equiv 0, S_{1} \equiv S_{2}$ and $S_{2} \equiv S_{3}$ as each of 
these would imply that some $[A_{i}]$ was zero.  Now, if $S_{1} = S_{3}$, 
we take the cycle representing their difference, and its complement.  
Since both cycles are connected, one of the two must satisfy the desired 
conditions.  If $S_{2}$ is zero mod 3, the same argument applies.  The 
hard case is when $S_{3} = 0$ integrally.

The cycle corresponding to $S_{3}$ is connected, but its complement, $A_{2} \cup A_{3}$ is not.  Since $S_{3} = 0$, we know that  $[A_{1}] = -([A_{4}] + [A_{5}])$.  Now, since $[A_{i}]$ is non-zero mod 3 for all $i$, $[A_{4}] \equiv [A_{5}]$ mod $3$ (otherwise $[A_{1}] \equiv 0$).  But now $[A_{1}] \equiv -2[A_{4}] \equiv [A_{4}]$ mod $3$.

Thus $[A_{1}] \equiv [A_{2}] \equiv [A_{4}] \equiv [A_{5}] \equiv \sum_{i\neq3} [A_{i}]$ mod $3$.  So, as $\sum [A_{i}] \equiv 0$ mod $3$, $[A_{3}] + [A_{5}] \equiv 0$ mod $3$.  Now both $A_{3} \cup A_{5}$ and its complement $A_{1} \cup A_{2} \cup A_{4}$ are connected, so one of the two is the desired cycle $A$.
\end{proof}

\begin{cor}Every embedding of $K_{7\alpha_{3n}'}$ contains an $n+1$ component link $L, Z_{1}, \ldots, Z_{n}$ with $lk(L,Z_{i}) = 3k$, $k \neq 0$. 
\label{mod3keys}
\end{cor}

\begin{proof}
The proof is essentially the same as the proof of Corollary \ref{mod2keys}.

We can think of $K_{7\alpha_{3n}'}$ as $\alpha_{3n}'$ copies of $K_{7}$. 
By Corollary \ref{starcor} we have a $3n+1$ component link $L,Z_{i}$ with $lk(L,Z_{i}) \neq 0$.  The $Z_{i}$ have at least four vertices each, so breaking them into groups of three, we may now apply construction from the proof of Theorem \ref{mod3thm} to each group individually to obtain the desired link.
\end{proof}

\Addresses\recd

\end{document}